\documentclass[11pt,a4paper]{article}  
\usepackage{amssymb}

\usepackage[dvips]{graphicx}
\usepackage[russian, english]{babel}
\begin{document}

\title{Asymptotic Properties of Hilbert Geometry}
\author{Alexandr A. BORISENKO and Eugeny A. OLIN
\footnote{The second author was partially supported by the Akhiezer Foundation}}
\date{23 April 2007}
\maketitle
\begin{center}
{\it Geometry Department, Mech.-Math. Faculty, Kharkov National
University, Pl. Svoboda 4, 310077-Kharkov, Ukraine.\\}
 E-mail: Alexander.A.Borisenko@univer.kharkov.ua, evolin@mail.ru
\end{center}

\begin{abstract}
We show that the spheres in Hilbert geometry have the same volume
growth entropy as those in the Lobachevsky space. We give the
asymptotic estimates for the ratio of the volume of metric ball to
the area of the metric sphere in Hilbert geometry. Derived
estimates agree with the well-known fact in the Lobachevsky space.
\end{abstract}

\textbf{Key words:} Hilbert geometry, Finsler geometry, balls, spheres,
volume, area, entropy.

\textbf{Mathematical Subject Classifications (2000):} 53C60, 58B20, 52A20

\section{Introduction}

Hilbert geometry is the generalization of the Klein model of the
Lobachevsky space. The absolute there is an arbitrary convex
hypersurface unlike an ellipsoid in the Lobachevsky space. Hilbert
geometries are simply connected, projectively flat, complete reversible Finsler spaces of constant negative flag curvature $-1$.

B. Colbois and P. Verovic proved in [12] that the balls in an
$(n+1)$-dimensional Hilbert geometry have the same volume growth
entropy as those in $\mathbb{H}^{n+1}$, namely $n$. We obtain the
analogous result for the spheres in Hilbert geometry.

\textbf{Theorem 1. }\textit{ Consider an $(n+1)$-dimensional
Hilbert geometry  associated with a bounded open convex domain $U
\subset \mathbb{R}^{n+1}$ whose boundary  is a $C^3$ hypersurface
with positive normal curvatures. Then we have}
$$\lim_{t\rightarrow \infty}\frac{\ln(\mathbf{Vol}(S_{t}^{n}))}{t}=n$$

It is known [4, 5, 6, 7] that in the Lobachevsky space
$\mathbb{H}^{n+1}$ of constant curvature $-1$ for a family of
metric balls $\{B_t^{n+1}\}_{t \in \mathbb{R}^+}$ the following
equality holds
$$\lim_{\rho\rightarrow\infty}\frac{\mathbf{Vol}(B_\rho^{n+1})}{\mathbf{Vol}(S_\rho^{n})} = \frac{1}{n}$$

Such a ratio in a more general case for $\lambda$- and $h$-convex
hypersurfaces in Hadamard manifolds was considered in [4, 6,
7] by A. A. Borisenko, V. Miquel, A. Reventos and E. Gallego.

Similar estimates in Finsler spaces were derived in [5] (see also [16]).

\textbf{Theorem [5].} \textit{Let $(M^{n+1},F)$ be an
$(n+1)$-dimensional Finsler-Hadamard manifold that satisfies the
following conditions:
\begin{enumerate}
\item Flag curvature satisfies the inequalities $-k_2^2\leqslant K
\leqslant -k_1^2$, $k_1,k_2>0$, \item $\mathbf{S}$-curvature
satisfies the inequalities $n\delta_1\leqslant S \leqslant
n\delta_2$ such that $\delta_i<k_i.$
\end{enumerate}}

\textit{Then for a family $\{B_r^{n+1}(p)\}_{r\geqslant0}$ we have
$$\frac{1}{n(k_2-\delta_2)}\leqslant \lim_{r\rightarrow\infty}\inf\frac{\mathbf{\mathbf{Vol}}( B_r^{n+1}(p))}{\mathbf{Area}(S_r^{n}(p))}\leqslant
\lim_{r\rightarrow\infty}\sup\frac{\mathbf{Vol}(
B_r^{n+1}(p))}{\mathbf{Area}(S_r^{n}(p))}\leqslant
\frac{1}{n(k_1-\delta_1)}.$$}

 Our goal is to prove analogous
result in Hilbert geometry for a family $\{B_t^{n+1}\}_{t \in
\mathbb{R}^+}$. Applying the theorem from [5] is the rather
difficult task because the $\mathbf{S}$-curvature in Hilbert
geometry is difficult to calculate.

As the result the following theorem is obtained.

\textbf{Theorem 2.} \textit{Consider an $(n+1)$-dimensional
Hilbert geometry  associated with a bounded open convex domain $U
\in \mathbb{R}^{n+1}$ whose boundary  is a $C^3$ hypersurface with
positive normal curvatures. Fix a point $o \in U$, we will
consider this point as the origin and the center of all the
considered balls. Denote by $\omega(u):\mathbb{S}^n\rightarrow
\mathbb{R}_{+}$ the radial function for $\partial U$, i. e. the
mapping $\omega(u)u$, $u \in \mathbb{S}^n$ is a parametrization of
$\partial U$, and by $\iota: \mathbb{R}^{n+1}\rightarrow \mathbb{S}^n$ the mapping such that $\iota(p) = \frac{u_p}{||u_p||}$, $u_p$ is the radius-vector of a point $p$.}

\textit{Denote by $K$ and $k$ the maximum and minimum normal
curvature of  $\partial U$, $c = \max_{u \in \mathbb{S}^n}
\frac{\omega(u)}{\omega(-u)}$, $\omega_0
= \min_{u \in \mathbb{S}^n} \omega(u)$, $\omega_1
= \max_{u \in \mathbb{S}^n} \omega(u)$. Then we have}

$$
\lim_{\rho\rightarrow\infty}\sup\frac{\mathbf{Vol}(B_\rho^{n+1})}{\mathbf{Vol}(S_\rho^{n})}
\leqslant \frac{1}{n} c^{\frac{n}{2}} \left(\frac{K}{k}\right)^{\frac{n}{2}} \frac{1}{(k \omega_0)^{\frac{n}{2}+1}}\frac{
\int_{\mathbb{S}^n}\omega(u)^\frac{n}{2}du}{
\int_{\partial U}\omega(\iota(p))^{-\frac{n}{2}}dp}$$

$$
\lim_{\rho\rightarrow\infty}\inf\frac{\mathbf{Vol}(B_\rho^{n+1})}{\mathbf{Vol}(S_\rho^{n})}
\geqslant \frac{1}{n} \frac{1}{c^{\frac{n}{2}}} \left(\frac{k}{K}\right)^{\frac{n}{2}} (k \omega_0)^{\frac{n}{2}} \frac{
\int_{\mathbb{S}^n}\omega(u)^\frac{n}{2}du}{
\int_{\partial U}\omega(\iota(p))^{-\frac{n}{2}}dp} $$

\textit{or, more simple expression}

$$
\lim_{\rho\rightarrow\infty}\sup\frac{\mathbf{Vol}(B_\rho^{n+1})}{\mathbf{Vol}(S_\rho^{n})}
\leqslant \frac{1}{n} \left(\frac{K}{k}\right)^{\frac{n}{2}} \left(\frac{\omega_1}{\omega_0}\right)^{n+1} \left(\frac{\omega_1}{k}\right)^{\frac{n}{2}} \frac{1}{k \omega _1}\frac{ \mathbf{Vol}_E(\mathbb{S}^n)}{\mathbf{Vol}_E(\partial U) }$$

$$
\lim_{\rho\rightarrow\infty}\inf\frac{\mathbf{Vol}(B_\rho^{n+1})}{\mathbf{Vol}(S_\rho^{n})}
\geqslant \frac{1}{n} \left(\frac{k}{K}\right)^{\frac{n}{2}} \left(\frac{\omega_0}{\omega_1}\right)^{\frac{n}{2}} \omega_0^n (k \omega_0)^{\frac{n}{2}} \frac{ \mathbf{Vol}_E(\mathbb{S}^n)}{\mathbf{Vol}_E(\partial U) } $$

\textit{If $U$ is a symmetric domain with respect to $o$ then we have}

$$
\lim_{\rho\rightarrow\infty}\sup\frac{\mathbf{Vol}(B_\rho^{n+1})}{\mathbf{Vol}(S_\rho^{n})}
\leqslant \frac{1}{n} c^{\frac{n}{2}} \left(\frac{K}{k}\right)^{\frac{n}{2}} \frac{\omega_1^n}{(k \omega_0)^{\frac{n}{2}+1}} \frac{ \mathbf{Vol}_E(\mathbb{S}^n)}{\mathbf{Vol}_E(\partial U) }$$

$$
\lim_{\rho\rightarrow\infty}\inf\frac{\mathbf{Vol}(B_\rho^{n+1})}{\mathbf{Vol}(S_\rho^{n})}
\geqslant \frac{1}{n} \frac{1}{c^{\frac{n}{2}}} \left(\frac{k}{K}\right)^{\frac{n}{2}} (k \omega_0)^{\frac{n}{2}} \omega_0^n\frac{ \mathbf{Vol}_E(\mathbb{S}^n)}{\mathbf{Vol}_E(\partial U) } $$

Notice that in this theorem the ratio of the
volume of the ball to the \textit{internal} volume of the sphere is considered,
unlike theorem [5], where the \textit{induced} volume is used.

\section{Preliminaries}
\subsection{Finsler geometry} In this section we recall some basic facts and theorems from
Finsler geometry that we need. See [16] for details.

 Let $M^n$ be an
$n$-dimensional connected $C^{\infty}$-manifold. Denote by
$TM^n=\bigsqcup_{x\in M^n}T_xM^n$ the tangent bundle of $M^n$,
where $T_xM^n$ is the tangent space at $x$. A \textit{Finsler
metric} on $M^n$ is a function $F:TM^n\rightarrow [0,\infty)$ with
the following properties:
\begin{enumerate}
\item $F\in C^{\infty}(TM^n\backslash\{0\})$;

\item $F$ is positively homogeneous of degree one, i. e. for any
pair $(x,y)\in TM^n$ and any $\lambda>0$, $F(x,\lambda y)=\lambda
F(x,y)$; \item For any pair $(x,y)\in TM^n$ the following bilinear
symmetric form $g_y:T_x M^n\times T_x M^n\rightarrow \mathbb{R}$
is positively definite,
$$\mathbf{g}_y(u,v):=\frac{1}{2}\frac{\partial^2}{\partial t \partial s}
\lbrack F^2(x,y+su+tv)\rbrack |_{s=t=0}$$
\end{enumerate}

The pair $(M^n,F)$ is called \textit{a Finsler manifold}.

If we denote by $$\mathbf{g}_{ij}(x,y) =
\frac{1}{2}\frac{\partial^2}{\partial y^i\partial
y^j}[F^2(x,y)],$$ then one can rewrite the form
$\mathbf{g}_y(u,v)$ as
$$\mathbf{g}_y(u,v)=\mathbf{g}_{ij}(x,y)u^iv^j$$

For any fixed vector field $Y$ defined on the subset $U \subset M^n$, $\mathbf{g}_Y(u,v)$ is a Riemannian
metric on $U$.

Given a Finsler metric $F$ on a manifold $M^n$. For a smooth curve
$c:[a,b]\rightarrow M^n$ the length is defined by the integral
$$L_F(c)=\int_a^b F(c(t),\dot{c}(t))dt = \int_a^b \sqrt{\mathbf{g}_{\dot{c}(t)}(\dot{c}(t),\dot{c}(t))}dt.$$

Let $\{e_i\}_{i=1}^n$ be an arbitrary basis for $T_xM^n$ and
$\{\theta^i\}_{i=1}^n$ the dual basis for $T_x^*M^n$. Consider the
set $B_{F(x)}^n=\left\{ (y^i)\in\mathbb{R}^n : F(x, y^ie_i)<1
\right \} \subset T_xM^n$. Denote by $\mathbf{Vol}_E(A)$ the
Euclidean volume of $A$. Then define the form
$$
dV_F=\sigma_F(x)\theta^1\wedge...\wedge\theta^n,
$$
here
\begin{equation}
\sigma_F(x):=\frac{\mathbf{Vol}_E(\mathbb{B}^n)}{\mathbf{Vol}_E(B_{F(x)}^n)}.
\end{equation} and $\mathbb{B}^n$ is the unit ball in
$\mathbb{R}^n$.

The volume form $dV_F$ determines a regular measure
$\mathbf{Vol}_F=\int dV_F$ and is called the
\textit{Busemann-Hausdorff volume form}.

For any Riemannian metric $g(u,v) = \mathbf{g}_{ij}(x)u^iv^j$ the
Busemann-Hausdorff volume form is the standard Riemannian volume
form
$$dV_g=\sqrt{\det(\mathbf{g}_{ij})}\theta^1\wedge...\wedge\theta^n.$$

It was proved in [9] that the  Busemann-Hausdorff measure for
reversible metric coincides with the $n$-dimensional outer
Hausdorff measure. Recall that the $n$\textit{-dimensional outer
Hausdorff measure of a set} $A$ is defined by
$$\nu_n=\lim_{r\rightarrow 0} \nu_{n,r},$$
$$\nu_{n,r} = \mathbf{Vol}_E(\mathbb{B}^n) \inf \left(\sum_{i} \rho_i^n : 2\rho_i<r, A \subseteq \bigcup_i B[x_i,\rho_i], x_i\in A\right)$$

It should be noticed here that if we calculate the Hausdorff
measure for the submanifold in a Finsler manifold with the
symmetric metric then we will obtain the \textit{internal} volume
on submanifold in the metric induced from the ambient space. But
unfortunately using of this volume implies certain difficulties.
In our case when we consider a sphere as the submanifold the
following claim does not hold

$$
\mathbf{Vol}(B_r^n)=\int_0^r\mathbf{Vol}(S_t^{n-1})dt,
$$
if we use the internal volume. For details, see [16].

\subsection{Hilbert geometry}

Consider a bounded open convex domain $U \subset \mathbb{R}^{n+1}$
whose boundary is a $C^3$ hypersurface with positive normal
curvatures in $\mathbb{R}^n$  equipped with a Euclidean norm $\|\cdot\|$.

For given two distinct points $p$ and $q$ in $U$, let $p_1$ and
$q_1$ be the corresponding intersection point of the half line
$p+\mathbb{R}_{-}(q-p)$ and $p+\mathbb{R}_{+}(q-p)$ with $\partial
U$ (Fig. 1).

\begin{figure}[h]
\begin{center}
\scalebox{0.3}[0.3]{
\includegraphics{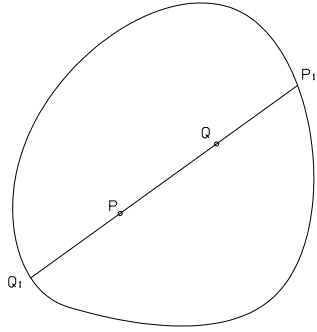}}
\caption{Hilbert metric} \end{center}
\end{figure}

Then consider the following distance function.

\begin{equation}d_U(p,q)=\frac{1}{2} \ln \frac{\|
q-q_1\|}{\|q-p_1\|}\times \frac{\|p-p_1\|}{\|p-q_1\|}
\end{equation}
$$d_U(p,p)=0$$
The obtained metric space $(U,d_U)$ is called Hilbert geometry
 and is a complete noncompact geodesic metric space with the $\mathbb{R}^n$-topology and in which the affine open
 segments joining two points are geodesics [10].

The distance function is associated in a natural way with the Finsler metric $F_U$
on $U$. For a point $p\in U$ and a tangent vector $v \in
T_pU=\mathbb{R}^n$
\begin{equation}
F_U(p,v)=\frac{1}{2}\|v\|\left(\frac{1}{\|p-p_{-}\|}+\frac{1}{\|p-p_{+}\|}\right)
\end{equation} where $p_{-}$ and $p_{+}$ is the intersection
point of the half-lines $p+\mathbb{R}_{-}v$ and
$p+\mathbb{R}_{+}v$ with $\partial U$.

Then $d_U(p,q)=\inf \int_IF_U(c(t),\dot{c}(t))dt$ when $c(t)$
ranges over all smooth curves joining $p$ to $q$.

In is known (see for example [16]) that Hilbert metrics are the
metrics of constant flag curvature $-1$.

 When $U=B^n_r$ then we obtain the Klein model of the $n$-dimensional Lobachevsky
 space $\mathbb{H}^n$ and the Finsler metric has the explicit
 expression
\begin{equation}
F_{B^n_r}(p,v)=\sqrt{\frac{\|v\|^2}{r-\|p\|^2}+\frac{<v,p>^2}{(r^2-\|p\|^2)^2}}
\end{equation}

It is proved in [10] that the balls of arbitrary radii are convex
sets in Hilbert geometry.

The asymptotic properties of Hilbert geometry have been obtained lately. All this properties mean that Hilbert geometry is "almost" Riemannian at infinity.
It is proved in [12] that Hilbert metric "tends" to riemannian
metric as follows.

\textbf{Theorem [12].} \textit{Let $\mathcal{C}\in \mathbb{R}^n$
be a bounded open convex domain whose boundary $\partial
\mathcal{C}$ is a hypersurface of class $C^3$ that is strictly
convex. For any $p\in \mathcal{C}$ let $\delta(p)>0$ be the
Euclidean distance from $p$ to $\partial \mathcal{C}$. Then there
exists a family $(\vec{l}_p)_{p\in \mathcal{C}}$ of linear
transformations in $\mathbb{R}^n$ such that
$$\lim_{\delta(p)\rightarrow0}\frac{F_C(p,v)}{\|\vec{l}_p(v)\|}=1$$
uniformly in $v\in \mathbb{R}^n\backslash\{0\}$}

This means that the unit sphere in the tangent space of given
Hilbert metric tends to ellipsoid in continuous topology as the tangent point goes to the absolute.

\section{Calculating the volume growth entropy of spheres.}

In this section we will prove that for an $(n+1)$-dimensional
Hilbert geometry
$$\lim_{t\rightarrow
\infty}\frac{\ln(\mathbf{Vol}(S_{t}^{n}))}{t}=n,$$ as it is in
$\mathbb{H}^{n+1}$.

Consider a bounded open convex domain $U \subset \mathbb{R}^{n+1}$
whose boundary is a $C^3$ hypersurface with positive normal
curvatures in $\mathbb{R}^n$.

Fix a point $o \in U$, we will consider this point as the origin
and the center of all the considered balls. Denote by
$\omega(u):\mathbb{S}^n\rightarrow \mathbb{R}_{+}$ the radial
function for $\partial U$, i. e. the mapping $\omega(u)u$, $u \in
\mathbb{S}^n$ is a parametrization of $\partial U$. Let
$B_r^{n+1}(o)$ be the metric ball of radius $r$ centered at a
point $o$, $S_r^{n}(o)=\partial B_r^{n+1}(o)$ be the metric
sphere.

We will use the following lemma that shows the order of growth
of the Hilbert distance from the sphere to $\partial U$ in terms of the
Euclidean distance. We also estimate the deviation of the tangent
and normal vectors to sphere from those to $\partial U$.

\textbf{Lemma 1.} \textit{Let $\omega(u)u:\mathbb{S}^n\rightarrow
\mathbb{R}_{+}$ be the parametrization of $\partial U$,
$\rho_{t}(u):\mathbb{S}^n\rightarrow \mathbb{R}_{+}$ -- the
parametrization of the sphere of radius $t$.}

\textit{Then, as $t\rightarrow\infty$:}
\begin{enumerate}
\item
$\omega(u)-\rho_t(u)=\Delta(u)e^{-2t}+\bar{o}(e^{-2t});$
$$\Delta(u) = \omega(u)\left(\frac{\omega(u)}{\omega(-u)} +
1 \right) $$
\item
$\omega_i'(u)-\rho_{t,i}'(u)=\Delta_i(u)e^{-2t}+\bar{o}(e^{-2t});$
$$\Delta_i(u) = \left[\omega'_i(u)\left(2\frac{\omega(u)}{\omega(-u)}+
1
\right)+\left(\frac{\omega(u)}{\omega(-u)}\right)^2\omega'_i(-u)\right] $$

\item
$\omega''_{ij}(u)-\rho_{t,ij}''(u)=\Delta_{ij}(u)e^{-2t}+\bar{o}(e^{-2t}),$
\end{enumerate}
$$\omega(u)^3\Delta_{ij}(u)= \omega(u)^2[2\omega_i'(-u)\omega_j'(-u)-\omega(-u)\omega_{ij}''(-u)] + \omega(-u)^2[2\omega_i'(u)\omega_j'(u)+\omega(-u)\omega_{ij}''(u)]+  $$
$$+2\omega(-u)\omega(u)[\omega_j'(-u)\omega_i'(u)+\omega_i'(-u)\omega_j'(u)] $$

P r o o f $ $ o f $ $ l e m m a 1. We are going to obtain the explicit
expression for  $\rho_{t}(u)$. Let $q=0$ be the center of the
sphere, $p$ be a point on the sphere. Using formula (3), we obtain
the equation on the function  $\rho_{t}(u)$
$$\frac{1}{2} \ln \left[\frac{\omega(u)}{\omega(-u)}\times
\frac{\omega(-u)+\rho_{t}(u)}{\omega(u)-\rho_{t}(u)}\right]=t$$

By the direct computation we have

$$\rho_{t}(u)=\frac{\omega(-u)\omega(u)(e^{2t}-1)}{\omega(u)+\omega(-u)e^{2t}}$$

\begin{enumerate}

\item Consider the difference
$$\omega(u)-\rho_{t}(u)=\omega(u)-\frac{\omega(-u)\omega(u)(e^{2t}-1)}{\omega(u)+\omega(-u)e^{2t}}=$$
$$=\frac{\omega^2(u)+\omega(-u)\omega(u)}{\omega(u)+\omega(-u)e^{2t}}=\omega(u)\left(\frac{\omega(u)}{\omega(-u)}
+ 1 \right)e^{-2t}+\bar{o}(e^{-2t}), t\rightarrow\infty$$

\item  We obtain analogously
$$\omega'(u)-\rho_{t,i}'(u)=$$
$$=\frac{\omega'_i(u)\omega(-u)^2e^{2t}+2e^{2t}\omega(u)\omega(-u)\omega'_i(u)+\omega(u)^2(\omega'_i(u)+\omega'_i(-u)(e^{2t}-1))}{(\omega(u)+\omega(-u)e^{2t})^2}=$$
$$=\left[\omega'_i(u)\left(2\frac{\omega(u)}{\omega(-u)}+ 1 \right)+\left(\frac{\omega(u)}{\omega(-u)}\right)^2\omega'_i(-u)\right] e^{-2t}+\bar{o}(e^{-2t}), t\rightarrow\infty$$

\item  It can be proved in the same manner.  $\square$
\end{enumerate}

 Denote by $k$ and $K$ the minimum and maximum Euclidean normal curvatures of $\partial U$.

 We also use the notations
$\omega_0=\min_{u \in \mathbb{S}^n}\omega(u)$, $\omega_1=\max_{u \in \mathbb{S}^n}\omega(u)$.

The following lemma gives the estimates on the angle between the
radial and normal directions at the points from $\partial U$.

\textbf{Lemma 2.} \textit{For a given point $m=\omega(u_m)u_m \in
\partial U$ denote by $N(m)$ the normal vector at $m$. Then}
$$\cos\angle(u_m,N(m))\geqslant \frac{\omega_0}{R}$$

P r o o f $ $ o f $ $ l e m m a 2. This lemma follows from the  more
general theorem.

\textbf{Theorem [4, 6, 7].} \textit{Let $N$ be a hypersurface in a
Riemannian manifold $M$. Consider $N$ as defined by the the equation
$t = \rho(\theta)$ of class $C^2$, where $\rho(\theta)$ is the distance to a point $o$.
$N$ can be seen as the 0-level set of the function $F=t-\rho$. For
given point $P\in N$ we consider all the vectors to be attached at
$P$. Denote by $Y = \frac{grad_N\rho}{\|grad_N\rho\|}$. Let $x$ be
a unit vector in the plane spanned on $y$ and the radial direction
that is orthogonal to the radial direction . Let $\varphi$ be
the angle between the normal direction and the radial direction at
the point $P \in N$.}

\textit{If $k_\mathbf{n}$ is the normal curvature at $P$ in
the direction given by $Y$, $\mu_\mathbf{n}$ is the normal
curvature in the direction of $x$ of the sphere centered at $o$ of
radius $\rho$ and $\frac{d\varphi}{ds}$ is the derivative of
$\varphi$ with respect to the arc parameter of the integral curve
of $Y$ by $P$ then}
$$k_\mathbf{n} = \mu_\mathbf{n}\cos\varphi+\frac{d\varphi}{ds}$$

Now we can prove lemma 2.

Consider any integral curve $\gamma$ of $\frac{y}{\|y\|}$. Since the
angle $\varphi$ takes its value in the interval $[0,\pi/2]$ then
there is a supremum $\varphi_0$ of it. If at some point
$\gamma(s_0)$ the value $\varphi_0$ is achieved then we have at
this point $\varphi'=0$ and
$$\cos\varphi=\frac{k_\mathbf{n}}{\mu_\mathbf{n}}$$

The minimum possible value of $k_\mathbf{n}$ is equal to $k =
\frac{1}{R}$, and the maximum possible value of $\mu_\mathbf{n}$
is equal to $\frac{1}{\omega_0}$. Hence we have

$$\cos\varphi=\frac{k_\mathbf{n}}{\mu_\mathbf{n}}\geqslant\frac{\omega_0}{R}$$

 And lemma 2 follows.
 $\square$

P r o o f $ $ o f $ $ t h e o r e m 1.  Now we are going to estimate the
volume of a sphere $S_{t}^{n}$ in Hilbert geometry. The idea of
proof is to obtain the Hausdorff measure of this sphere. It
follows from the reversibility of Hilbert metrics that the
Hausdorff measure coincides with the Finslerian
Busemann-Hausdorff volume [9].

Fix the point $p$ on the sphere $S_r^n$. Since the spheres are convex we can choose the
vector $u\in \mathbb{S}^n$ such that $p = \rho_t(u)$. More generally, for a given origin $o\in \mathbb{R}^{n+1}$ denote by $u_p$ the corresponding radius vector and consider the function  $\iota: \mathbb{R}^{n+1}\rightarrow \mathbb{S}^n$ such that $\iota(p) = \frac{u_p}{||u_p||}$. Then we can write that $p = \rho_t(\iota(u))$.

Denote by $m$ the point $\omega(\iota(p))\iota(p) \in \partial U$. Consider the vector $v_m$ which is tangent to $\partial U$ at $m$, the vector $n_m$ which is orthogonal to $v_m$ with respect to the Euclidean inner product such that the point $o$ belongs to the plane $\mathcal{P}$ spanned on $v_m$ and $n_m$. Let $k_m$ be the curvature of the section of $\partial U$ by $\mathcal{P}$ at $m$. Consider the special coordinate system in the plane $\mathcal{P}$: let the axe $z$ be directed as $n_m$, and the axe $x$ be directed as $v_m$.
Then in this special coordinate system the section of $\partial U$ can be locally
expressed as

$$z(x) = \frac{1}{2} k_m x^2 + \bar{o}(x^2), x\rightarrow 0$$
Later on we will work with this section.

 Draw the secant of the sphere that is
parallel to the tangent vector at $p$ (Fig. 2).

\begin{figure}[h]
\begin{center}
\scalebox{0.65}[0.65]{
\includegraphics{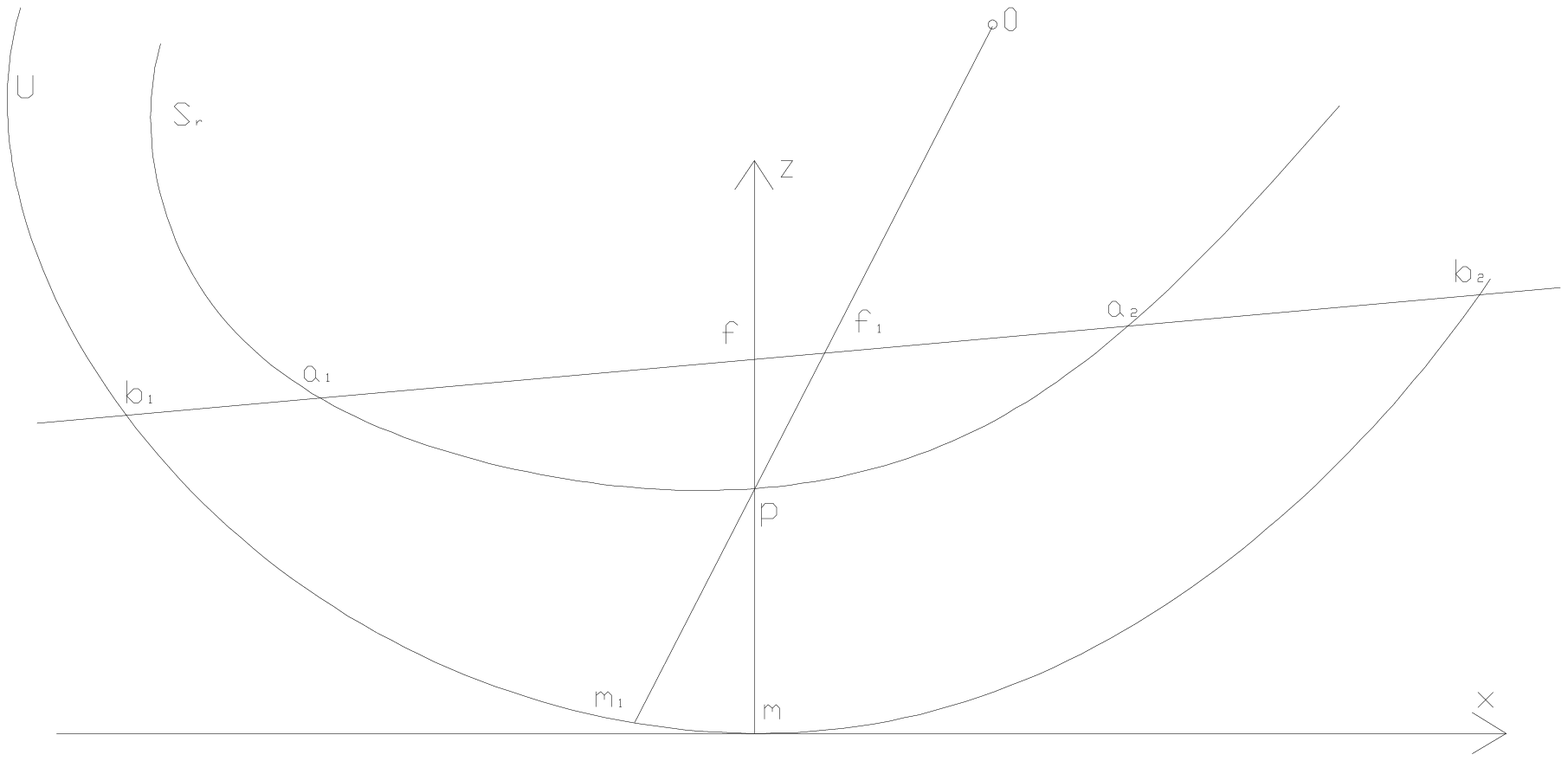}}
\caption{} \end{center}
\end{figure}

Put $d = ||a_1-a_2||$, $\delta(p) = ||m-p||$, $\delta_1 = ||b_1-f||$,
$\delta_2 = ||f-b_2||$, $h = ||p-f||$.

Let us estimate the function $\delta(p)$. From lemma 1 we have
\begin{equation} \delta(p) \leqslant \omega(u)\left(\frac{\omega(u)}{\omega(-u)} +
1 \right)e^{-2t}+\bar{o}(e^{-2t}), t\rightarrow\infty \end{equation}

From the triangle $pm_1m$ we have that $\delta(p) \thickapprox \cos\angle(u_m,v_m)
(\omega(u)-\rho_t(u))$. Finally, using lemma 2 we obtain $\delta(p) \geqslant
\frac{\omega_0}{R}(\omega(u)-\rho_t(u))$

Consequently, \begin{equation} \delta(p) \geqslant \frac{\omega_0}{R}\omega(u)
\left(\frac{\omega(u)}{\omega(-u)} + 1
\right)e^{-2t}+\bar{o}(e^{-2t}), t\rightarrow \infty\end{equation}

 Then we estimate $\delta_1$ and $\delta_2$. Let
$$z = a(p)x + \delta(p) + h $$ be the equation of the secant in the special coordinates system.
 We will think at once that  $h$
decreases faster than $\delta(p)$. Thus in the further
computations we will neglect $h$.

We find the intersection points of this line with the boundary
$\partial U$.
 From the expression for the boundary we have $$a(p)x+\delta(p) =
 \frac{1}{2}k_mx^2.$$

 Thus

$$x_{1,2} = \frac{a(p)\pm \sqrt{a(p)^2 + 2k_m \delta(p)}}{k_m}$$

It follows from lemma 1 that $a(p) = a_0(u)e^{-2t}, t
\rightarrow\infty $, for some function $a_0(u)$ and consequently
$a(p) = \underline{\mathrm{O}}(\delta(p)), \delta(p) \rightarrow 0$. Therefore we have

$$x_{1,2} = \pm  \sqrt{\frac{2\delta(p)}{k_m}}+\bar{o}(\sqrt{\delta(p)}),\delta(p)\rightarrow 0$$

$$z_{1,2} = \frac{1}{2} k_m x^2 + \bar{o}(x^2)|_{x=x_{1,2}} = \frac{1}{2}\delta(p)+\bar{o}(\delta(p)),\delta(p)\rightarrow0$$

and

$$\delta_1  = \sqrt{x_1^2+(z_1-\delta(p))^2} = \sqrt{\frac{2\delta(p)}{k_m}}+\bar{o}(\delta(p))= \delta_2,\delta(p)\rightarrow 0$$

Therefore the turning of the tangent as the point goes to the
boundary does not influence on the asymptotic behavior of
$\delta_i$.

Compute the Hilbert length of the segment $a_1a_2$. Denote it by
$d_U$. Then as $h \rightarrow 0$:
$$d_U\thickapprox\frac{1}{2}\ln\left[\frac{d+\delta_1}{\delta_2}\times\frac{d+\delta_2}{\delta_1}\right]
=
\frac{1}{2}\ln\left[\left(\frac{d}{\delta_2}+\frac{\delta_1}{\delta_2}\right)\times\left(\frac{d}{\delta_1}+\frac{\delta_2}{\delta_1}\right)\right]\thickapprox$$
$$\thickapprox \frac{1}{2}\left(\frac{d}{\delta_1}+\frac{d}{\delta_2}\right) \thickapprox \frac{d \sqrt{k_m}}{\sqrt{2(\delta(p)+h)}} + \bar{o}(\sqrt{1/\delta(p)}),\delta(p)\rightarrow 0$$

We are showing now that the limit of the ratio of $d_U$ to the
Finslerian length $\tilde{d}_U$ of the geodesic arc $a_1a_2$ is equal to
$1$ as the arc is subtended to a point. Specialize the
coordinate system on $\mathbb{R}^{n+1}$ so as $a_1=0$. Let
$w(t):[0,T]\rightarrow U$ be a parametrization of the arc. Then
the segment from the point $a_1=0$ to the point $a_2=w(t)$ can be
parameterized by $v(s) = \frac{s}{t}w(t) :[0,t]\rightarrow U $.
Calculate lengths of $v$ and $w$.

$$\tilde{d}_U = \int_0^t F_U(w(s),\dot{w}(s))ds$$
$$d_U = \int_0^t F_U(v(s),\dot{v}(s))ds = \int_0^t F_U\left(\frac{s}{t}w(t),\frac{1}{t}w(t)\right)ds.$$

From the intermediate-value theorem for integrals we have
$$\tilde{d}_U = \int_0^t F_U(w(s),\dot{w}(s))ds = t F_U(w(s_0),\dot{w}(s_0)), s_0\in [0,t]$$
$$d_U = \int_0^t F_U\left(\frac{s}{t}w(t),\frac{1}{t}w(t)\right)ds = t F_U\left(\frac{s_1}{t}w(t),\frac{1}{t}w(t)\right), s_1\in [0,t]$$

Now we are subtending the arc to a point, i. e. let
$t\rightarrow0$. Then $s_0,s_1\rightarrow0$, and we obtain:

$$\frac{\tilde{d}_U}{d_U} = \frac{t F_U(w(s_0),\dot{w}(s_0))}{t F_U\left(\frac{s_1}{t}w(t),\frac{1}{t}w(t)\right)}\longrightarrow\frac{F_U(0,\dot{w}(0))}{F_U\left(0,\dot{w}(0)\right)} = 1$$

And the statement is proved.

Now our goal is to calculate the Hausdorff measure of the sphere $S_r^n$. Denote by $\delta_0(r)$ the Hausdorff distance
from the points of the sphere to the absolute $\partial U$. Consider a covering $\{B_i\}$ of the sphere $S_r^n$ by balls of diameters $\tilde{d}_i$ centered at points $p_i \in S_r^n$. Denote by $k_i$ the normal curvature of $\partial U$ that corresponds to the $i$-th sphere from the covering (as above). As we saw, we can replace $\tilde{d}_i$ by the lengths of the corresponding chords $d_i$ of the sphere $S_r^n$.

Then the Hausdorff measure, and, consequently, the Finslerian Busemann-Hausdorff measure is given by
$$\mathbf{Vol}(S_t^n)= \mathbf{Vol}_E(\mathbb{B}^n) \inf_{d_{B_i}}\sum_i\left(\frac{d_i
\sqrt{k_i}}{\sqrt{2(\delta(p_i)+h)}}\right)^n +
\bar{o}(\sqrt{1/\delta_0(t)^n}),\delta_0(t)\rightarrow 0,$$

where infimum is calculated over all  coverings of the sphere $S_r^n$.

Our metric sphere $S_r^n$ is sufficiently smooth, so we can proceed to the integral over $S_r^n$.

$$\mathbf{Vol}(S_t^n)= \mathbf{Vol}_E(\mathbb{B}^n) \inf_{d_{B_i}}\sum_i\left(\frac{d_i
\sqrt{k_i}}{\sqrt{2(\delta(p_i)+h)}}\right)^n +
\bar{o}(\sqrt{1/\delta_0(t)^n}) = $$
$$  = \inf_{d_{Bi}}\sum_i\left(\frac{
\sqrt{k_i}}{\sqrt{2(\delta(p_i)+h)}}\right)^n \mathbf{Vol}_E(B_i) +
\bar{o}(\sqrt{1/\delta_0(t)^n}) ,\delta_0(t)\rightarrow 0,$$

Denote by $dp$ the area element of $S_t^n$. Proceeding to integral and estimating leads to
$$\mathbf{Vol}(S_t^n)\geqslant k^{\frac{n}{2}}\int_{S_t^n}\left(\frac{1}{2\delta(p)}\right)^{\frac{n}{2}}dp +
\bar{o}(\sqrt{1/\delta_0(t)^n}),\delta_0(t)\rightarrow 0 $$
$$\mathbf{Vol}(S_t^n)\leqslant K^{\frac{n}{2}}\int_{S_t^n}\left(\frac{1}{2\delta(p)}\right)^{\frac{n}{2}}dp +
\bar{o}(\sqrt{1/\delta_0(t)^n}),\delta_0(t)\rightarrow 0 $$

 Using the explicit estimates (5), (6) for $\delta(p)$ as the result we
have

\begin{equation}
\mathbf{Vol}(S_{t}^n) \geqslant
 k^{\frac{n}{2}}\int_{S_t^n}\left(2\omega(\iota(p))\left(\frac{\omega(\iota(p))}{\omega(-\iota(p))}
+ 1 \right)\right)^{-\frac{n}{2}}dp \cdot e^{n t}+\bar{o}(e^{n t}), t\rightarrow \infty
\end{equation}
\begin{equation}
\mathbf{Vol}(S_{t}^n)\leqslant
 K^{\frac{n}{2}}\int_{S_t^n}\left(2
\frac{\omega_0}{R}\omega(\iota(p)) \left(\frac{\omega(\iota(p))}{\omega(-\iota(p))} + 1
\right)\right)^{-\frac{n}{2}}dp\cdot e^{n t} + \bar{o}(e^{n t}) ,
t\rightarrow \infty
\end{equation}

And theorem 1 follows. $\blacksquare$

\section{Estimation of the ratio of the volume of the ball to the volume of the sphere.}

Here we will find the asymptotic behavior of the volume of the
metric ball $B_\rho^{n+1}$ in Hilbert geometry. We will use the
method introduced in [13] in which some necessary estimates were
improved.

The volume of a metric ball is given by the integral
$$\mathbf{Vol}(B_\rho^{n+1})=\int_{B_\rho^{n+1}}\sigma(p)dp$$
Here $\sigma(p)$ is the Busemann-Hausdorff volume form. And the
volume estimating problem is reduced to the estimating of the
volume form. Recall (1) that
$$\sigma(p) :=
\sigma_{F_U}(p)=\frac{\mathbf{Vol}_E(\mathbb{B}^n)}{\mathbf{Vol}_E(B_{F_U(p)}^n)}.$$
Thus we have to estimate the volume of the unit sphere in the
tangent space at the point $p \in U$.

We will use the following simple lemma.

\textbf{Lemma 3.} \textit{There exists a value $\rho_0$ such that
for any points $p \in U$ in the neighborhood
$d(p,\partial U)\leqslant\rho_0$ there exist a unique point $\pi(p) \in \partial U$:
$d(p,\pi(p))=d(p,\partial U)$ }

Put $m = \pi(p) \in
\partial U$. Denote by $k$ and $K$ the minimum and maximum Euclidean normal curvatures of $\partial U$.
Then at any point $m \in \partial U$ the tangent sphere of radius
$R := \frac{1}{k}$ contains $U$, the tangent sphere of radius $r :=
\frac{1}{K}$ is contained in $U$ [2]. On two tangent spheres of the radii $r$ and $R$ at this point we construct corresponding Klein metrics $F_r$
and $F_R$.
 We can give the explicit expressions (4) for them.

Then the following inequalities hold
\begin{equation}\mathbf{Vol}_E\left(B_{F_r(p)}^{n+1}\right)\leqslant
\mathbf{Vol}_E\left(B_{F_U(p)}^{n+1}\right)\leqslant
\mathbf{Vol}_E\left(B_{F_R(p)}^{n+1}\right)\end{equation}

As it was shown in [13]:
$$
\mathbf{Vol}_E\left(B_{F_R(p)}^{n+1}\right)=\mathbf{Vol}_E(\mathbb{B}^{n+1})
R^{n+1}
\left\{1-\left(1-\frac{d(p,m)}{R}\right)^2\right\}^{\frac{n+2}{2}}
$$
$$
\mathbf{Vol}_E\left(B_{F_r(p)}^{n+1}\right)=\mathbf{Vol}_E(\mathbb{B}^{n+1})
r^{n+1}
\left\{1-\left(1-\frac{d(p,m)}{r}\right)^2\right\}^{\frac{n+2}{2}}
$$

Thus, we have
\begin{equation}
\frac{1}{R^{n+1}
\left\{1-\left(1-\frac{d(p,m)}{R}\right)^2\right\}^{\frac{n+2}{2}}}  \leqslant \sigma(p)  \leqslant \frac{1}{r^{n+1}
\left\{1-\left(1-\frac{d(p,m)}{r}\right)^2\right\}^{\frac{n+2}{2}}}
\end{equation}

Consider the mapping
$$\Phi(u,s)=\tanh(s)\omega(u)u:\mathbb{S}^n \times \mathbb{R}\longrightarrow U$$
It was shown in [13] that the mapping $\Phi(u,s)$ satisfies the
following properties
\begin{enumerate}
\item $\Phi(\mathbb{S}^n,[0,\rho-c]) \subseteq
B_{\rho}^{n+1}\subseteq \Phi(\mathbb{S}^n,[0,\rho+1])$ where $c =
\sup_{u \in \mathbb{S}^n} \frac{\omega(u)}{\omega(-u)}$

Hence, $$\mathbf{Vol}(\Phi(\mathbb{S}^n,[0,\rho-c])) \leqslant
\mathbf{Vol}(B_\rho^{n+1})\leqslant
\mathbf{Vol}(\Phi(\mathbb{S}^n,[0,\rho+1]))
$$

\item
$|\mathbf{Jac}(\Phi(u,s))|=\omega(u)^{n+1}\tanh^{n}(s)(1-\tanh^2(s))$
\end{enumerate}

We improve the first property.

Fix $d>0$. Consider the difference
$$\rho_t(u) - \omega(u)\tanh(t+d) = \omega(u)\left(1-\frac{\omega(u)+\omega(-u)}{e^{2t}\omega(-u)+\omega(u)}-\tanh(t+d)\right)=$$
$$=\omega(u)\left( \frac{2}{e^{2(t+d)}+1} -
\frac{\omega(u)+\omega(-u)}{e^{2t}\omega(-u)+\omega(u)}\right) =$$
$$=\rho_t(u) - \omega(u)\tanh(t+d)=
\omega(u)e^{-2t}\left( 2e^{-2d} - 1 -
\frac{\omega(u)}{\omega(-u)}\right) +
\bar{o}(e^{-2t}),t\rightarrow\infty$$

Thus $B_{\rho}^{n+1}\subseteq \Phi(\mathbb{S}^n,[0,\rho+d])$ for
sufficiently large  $\rho$ if
$$2e^{-2d} - 1 -
\frac{\omega(u)}{\omega(-u)} \leqslant 0$$

$$d \geqslant -\frac{1}{2}\ln\left[\frac{1}{2}\left(1+\frac{1}{c}\right)\right] := d_1$$
and  $B_{\rho}^{n+1}\supseteq \Phi(\mathbb{S}^n,[0,\rho+d])$ for
sufficiently large finite  $\rho$ if

$$d \leqslant -\frac{1}{2}\ln\left[\frac{1}{2}\left(1+c\right)\right]:=d_2$$

Fix the values $d_2 $ and $d_1$ and choose sufficiently large $\rho_0$.

Then
\begin{equation}
\mathbf{Vol}(\Phi(\mathbb{S}^n,[0,\rho+d_2])) \leqslant
\mathbf{Vol}(B_\rho^{n+1})\leqslant
\mathbf{Vol}(\Phi(\mathbb{S}^n,[0,\rho+d_1]))
\end{equation}

Notice that if the domain $U$ is centrally-symmetric then $d_1 = d_2 = 0$.
In the worst case when $c \rightarrow \infty$ we have $d_1
\rightarrow \ln \sqrt{2} \approx 0.347 <1 $. Inclusion
(11) is more precise  than it was obtained in [13]. It will be
essentially used in the proof of theorem  2.

The volume of the set $\Phi(\mathbb{S}^n,[\rho_0,\rho])$ is given
by.
$$\mathbf{Vol}(\Phi(\mathbb{S}^n,[\rho_0,\rho]))=\mathbf{Vol}_E(\mathbb{B}^n)\int_{\mathbb{S}^n}\int_{\rho_0}^{\rho}\sigma(\Phi(u,s))|\mathbf{Jac}(\Phi(u,s))|dsdu$$

It is known [13] that
$$|Jac(\Phi(u,s))|=\omega(u)^{n+1}\tanh^{n}(s)(1-\tanh^2(s))=\omega(u)^{n+1} \frac{4e^{2s}\left(\frac{e^{2s}-1}{e^{2s}+1}\right)^{n+1}}{e^{4s}-1}$$

And, using the estimates (10) we obtain
$$\int_{\mathbb{S}^n}\int_{\rho_0}^{\rho}\frac{4 \omega(u)^{n+1} \frac{e^{2s}\left(\frac{e^{2s}-1}{e^{2s}+1}\right)^{n+1}}{e^{4s}-1}}{R^{n+1}\left(1-\left(1-\frac{d(\Phi(u,s),\partial U)}{R}\right)^2\right)^{\frac{n+2}{2}}}dsdu \leqslant \mathbf{Vol}(\Phi(\mathbb{S}^n,[\rho_0,\rho]))$$
$$\mathbf{Vol}(\Phi(\mathbb{S}^n,[\rho_0,\rho]))  \leqslant  \int_{\mathbb{S}^n}\int_{\rho_0}^{\rho}\frac{4 \omega(u)^{n+1} \frac{e^{2s}\left(\frac{e^{2s}-1}{e^{2s}+1}\right)^{n+1}}{e^{4s}-1}}{r^{n+1}\left(1-\left(1-\frac{d(\Phi(u,s),\partial U)}{r}\right)^2\right)^{\frac{n+2}{2}}}dsdu$$

Out next task is to find the asymptotic behavior of the integral
$$\int_0^r\frac{\frac{4e^{2s}\left(\frac{e^{2s}-1}{e^{2s}+1}\right)^{n+1}}{e^{4s}-1}}{(1-(1-Ce^{-2s})^2)^{\frac{n+2}{2}}}ds$$
After the changing of the variable $y=e^{-2s}$, we obtain the
integral
$$\int_{e^{-2r}}^{1} \frac{-8\frac{y^{-2}\left(\frac{y^{-1}-1}{y^{-1}+1}\right)^{n+1}}{y^{-2}-1}}{(1-(1-Cy)^2)^{\frac{n+2}{2}}} dy = \int_{e^{-2r}}^{1} \frac{8(1-y)^{n+1}}{(1+y)^{n+1}(y^2-1)(Cy(2-Cy))^{\frac{n+2}{2}}}dy= $$
$$=\int_{e^{-2r}}^{1}
\frac{1}{C^{\frac{n+2}{2}}y^{\frac{n+2}{2}}2^{\frac{n}{2}-2}} \cdot
\frac{(1-y)^{-\frac{n}{2}}2^{\frac{n+2}{2}}}{(1+y)^{n+1}(y^2-1)(2-Cy)^{\frac{n+2}{2}}}dy$$
Notice that $$\lim_{y\rightarrow
0}\left[\frac{(1-y)^{-\frac{n}{2}}2^{\frac{n+2}{2}}}{(1+y)^{n+1}(y^2-1)(2-Cy)^{\frac{n+2}{2}}}\right]=-1$$

Taking this into account and making the inverse change of variable
we get
\begin{equation}
\int_0^r\frac{\frac{4e^{2s}\left(\frac{e^{2s}-1}{e^{2s}+1}\right)^{n+1}}{e^{4s}-1}}{(1-(1-Ce^{-2s})^2)^{\frac{n+2}{2}}}ds
= \frac{1}{n C^{\frac{n+2}{2}}2^{\frac{n-2}{2}}}e^{nr} +
\bar{o}(e^{nr}),r\rightarrow\infty
\end{equation}

The expression for $\mathbf{Vol}_E\left(B_{F_R(p)}^{n+1}\right)$
includes the quantity $d(p,m) = d(p,\partial U)$. Thus we need the estimates of
$d(p,m)$ for the point $p = \Phi(u,s)$. So,
$$
d(\Phi(u,s),\omega(u)u) = \omega(u)-\tanh(s)\omega(u) = \omega(u)
- \frac{e^{2s}-1}{e^{2s}+1}\omega(u) = \frac{2\omega(u)}{1+e^{2s}}
$$ Finally,
\begin{equation}
d(\Phi(u,s),\partial U) \leqslant
2\omega(u)e^{-2s}+\bar{o}(e^{-2s}) \end{equation}

On the other hand analogously as formula (6) we get

 \begin{equation} d(\Phi(u,s),\partial U) \geqslant
2{\frac{\omega_0}{R}}\omega(u)e^{-2s}+\bar{o}(e^{-2s})
\end{equation}

Using (12), (13), (14), one can compute that
\begin{equation}
\frac{1}{n} \mathbf{C}_1 e^{n \rho } + \bar{o}(e^{n \rho })
\leqslant \mathbf{Vol}(\Phi(\mathbb{S}^n,[\rho_0,\rho])) \leqslant
\frac{1}{n} \mathbf{C}_2 e^{n \rho } + \bar{o}(e^{n \rho
}),\rho\rightarrow\infty
\end{equation}
$$\mathbf{C}_1 = \frac{1}{2^n} \int_{\mathbb{S}^n}\left(\frac{\omega(u)}{R}\right)^\frac{n}{2}du $$
$$\mathbf{C}_2 = \frac{1}{2^n}
\frac{R^{\frac{n+2}{2}}}{\omega_0^{\frac{n+2}{2}}}
\int_{\mathbb{S}^n}\left(\frac{\omega(u)}{r}\right)^\frac{n}{2}du
$$

And, taking into account (11), (15), we have
\begin{equation}
\frac{1}{n} \mathbf{C}_1 e^{nd_2} e^{n \rho } + \bar{o}(e^{n \rho
}) \leqslant \mathbf{Vol}(B_\rho^{n+1}) \leqslant \frac{1}{n}
\mathbf{C}_2 e^{n \rho }e^{nd_1} + \bar{o}(e^{n \rho
}),\rho\rightarrow\infty
\end{equation}

P r o o f $ $ o f $ $ t h e o r e m 2. It follows from (6), (7), (16) that:

$$
\lim_{\rho\rightarrow\infty}\sup\frac{\mathbf{Vol}(B_\rho^{n+1})}{\mathbf{Vol}(S_\rho^{n})}
\leqslant \frac{1}{n}\frac{1}{2^{n/2}}e^{nd_1}
\frac{R^{\frac{n+2}{2}}}{\omega_0^{\frac{n+2}{2}}} \frac{
\int_{\mathbb{S}^n}\left(\frac{\omega(u)}{r}\right)^\frac{n}{2}du}{
k^{\frac{n}{2}}\int_{\partial U}\left(\omega(\iota(p))\left(\frac{\omega(\iota(p))}{\omega(-\iota(p))}
+ 1 \right)\right)^{-\frac{n}{2}}dp} $$
$$\leqslant \frac{1}{n} c^{\frac{n}{2}} \left(\frac{K}{k}\right)^{\frac{n}{2}} \frac{1}{(k \omega_0)^{\frac{n}{2}+1}}\frac{
\int_{\mathbb{S}^n}\omega(u)^\frac{n}{2}du}{
\int_{\partial U}\omega(\iota(p))^{-\frac{n}{2}}dp}  $$
$$ \leqslant \frac{1}{n} c^{\frac{n}{2}} \left(\frac{K}{k}\right)^{\frac{n}{2}} \frac{\omega_1^n}{(k \omega_0)^{\frac{n}{2}+1}} \frac{ \mathbf{Vol}_E(\mathbb{S}^n)}{\mathbf{Vol}_E(\partial U) }$$

Note that $c \leqslant \frac{\omega_1}{\omega_0}$. Hence
$$ \lim_{\rho\rightarrow\infty}\sup\frac{\mathbf{Vol}(B_\rho^{n+1})}{\mathbf{Vol}(S_\rho^{n})} \leqslant \frac{1}{n} \left(\frac{K}{k}\right)^{\frac{n}{2}} \left(\frac{\omega_1}{\omega_0}\right)^{n+1} \left(\frac{\omega_1}{k}\right)^{\frac{n}{2}} \frac{1}{k \omega _1}\frac{ \mathbf{Vol}_E(\mathbb{S}^n)}{\mathbf{Vol}_E(\partial U) }$$

$$
\lim_{\rho\rightarrow\infty}\inf\frac{\mathbf{Vol}(B_\rho^{n+1})}{\mathbf{Vol}(S_\rho^{n})}
\geqslant \frac{1}{n}\frac{1}{2^{n/2}}e^{nd_2} \frac{
\int_{\mathbb{S}^n}\left(\frac{\omega(u)}{R}\right)^\frac{n}{2}du}{
K^{\frac{n}{2}}\int_{\partial U}\left(
\frac{\omega_0}{R}\omega(\iota(p)) \left(\frac{\omega(\iota(p))}{\omega(-\iota(p))} + 1
\right)\right)^{-\frac{n}{2}}dp}$$
$$ \geqslant \frac{1}{n} \frac{1}{c^{\frac{n}{2}}} \left(\frac{k}{K}\right)^{\frac{n}{2}} (k \omega_0)^{\frac{n}{2}} \frac{
\int_{\mathbb{S}^n}\omega(u)^\frac{n}{2}du}{
\int_{\partial U}\omega(\iota(p))^{-\frac{n}{2}}dp}$$

$$\geqslant \frac{1}{n} \frac{1}{c^{\frac{n}{2}}} \left(\frac{k}{K}\right)^{\frac{n}{2}} \omega_0^n (k \omega_0)^{\frac{n}{2}} \frac{ \mathbf{Vol}_E(\mathbb{S}^n)}{\mathbf{Vol}_E(\partial U) }\geqslant \frac{1}{n} \left(\frac{k}{K}\right)^{\frac{n}{2}} \left(\frac{\omega_0}{\omega_1}\right)^{\frac{n}{2}} \omega_0^n (k \omega_0)^{\frac{n}{2}} \frac{ \mathbf{Vol}_E(\mathbb{S}^n)}{\mathbf{Vol}_E(\partial U) }  $$

And the theorem follows. $\blacksquare$

E x a m p l e 1. Let $U = \mathbb{B}_{\rho}^{n+1}$. Then we get the
Klein model of the Lobachevsky space. Applying theorem 2 to this
space implies
$$\omega(u) = \frac{1}{k} = \frac{1}{K} = r = R = \omega_0 = \rho$$
$$c = 1$$
$$\int_{\partial U}du = \rho^n \mathbf{Vol}_E(\mathbb{S}^n)$$

Therefore we have obtained the well-known result

$$\lim_{\rho\rightarrow\infty}\frac{\mathbf{Vol}(B_\rho^{n+1})}{\mathbf{Vol}(S_\rho^{n})} = \frac{1}{n}$$

E x a m p l e 2. One  should not hope that for all  metrics of
negative curvature such result holds.

 Let $U$ be a open bounded strongly convex domain
in $\mathbb{R}^n$, $o = 0\in \mathbb{R}^n$. Given a point $x \in U$ and a direction $y \in
T_xU\backslash\{0\}\simeq U\backslash\{0\}$. The \textit{Funk
metric $F(x,y)$} is a Finsler metric that satisfies the following
condition
$$x+\frac{y}{F(x,y)}\in \partial U.$$

Then Hilbert metric is a symmetrized Funk metric

$$F_U(x,y)=\frac{1}{2}\left[F(x,y)+F(x,-y)\right]$$

 Funk metrics are  of constant negative curvature $-1/4$, but for such metrics
 [5]:

$$\lim_{r\rightarrow\infty}\frac{\mathbf{Vol}(
B_r^{n+1})}{\mathbf{Vol}(S_r^{n})}=\infty.$$

\end{document}